%
%
\input amstex
\documentstyle{tmna}
\NoBlackBoxes
\leftheadtext{Daciberg L. Gon\c calves and Ulrich Koschorke}
\rightheadtext{Coincidence theory of fibre-preserving maps}

\topmatter
\title Nielsen coincidence theory of fibre-preserving maps   and Dold's fixed point index\endtitle
\author Daciberg L. Gon\c calves $^1$  and Ulrich Koschorke $^2$ \endauthor
\address {\smc Daciberg L. Gon\c calves}\\
Dept. de Matem\'atica - IME - USP, 
Caixa Postal 66.281 - CEP 05311- 970,
 S\~ao Paulo - SP, BRASIL
\endaddress 
 \email dlgoncal\@ime.usp.br\endemail
\address {\smc Ulrich Koschorke}\\
Department of Mathematics, Universit\"at Siegen,
 Emmy-Noether-Campus,
D-57068 Siegen, DEUTSCHLAND
 \email koschorke\@mathematik.uni-siegen.de\endemail
\endaddress



\dedicatory Dedicated to Albrecht Dold on the occasion of his 80th birthday \enddedicatory

\subjclass Primary 55M20; Secondary 55R05, 55S35, 57R19, 57R90\endsubjclass

\keywords 
coincidence, fixed point,   map over $B$,  normal bordism, $\omega$-invariant,  Nielsen number, Reidemeister class, Dold's index, fibration \endkeywords

\abstract Let $M \to B$, $N \to B$ be fibrations and $f_1,f_2:M \to N$ be a pair of  fibre-preserving maps. Using normal bordism techniques we define an invariant  which is an obstruction to deforming the pair  $f_1,f_2$ over $B$ to a coincidence  free pair of maps. 
In the special case where the two fibrations are the same and one of the maps is the identity, a weak version of our 
$\omega$-invariant  turns out to equal Dold's fixed point index of fibre-preserving maps. The concepts of Reidemeister classes
and Nielsen coincidence classes  over  $B$ are developed. As an illustration we compute e.g. the minimal number of coincidence components for all homotopy classes of maps between  $S^1-$bundles over $S^1$ as well as their Nielsen and  Reidemeister numbers.   
 \endabstract

\thanks This work was  initiated during the visit of the second  author to the Institute de Matem\'atica e Estat\'istica-USP  September 23- October 21, 2007. The visit was supported by the international cooperation program DAAD-Capes. 
\endthanks
\endtopmatter

\document

\head 1. Introduction and outline of  results \endhead 

Throughout this paper we consider the following situation:



$$\matrix
& F_M & & & & F_N\\
&\bigcap & & & & \bigcap\\
f,\,f_1,\,f_2,\,f^\prime,\,f_1^\prime,\,f_2^\prime,\,\ldots\colon & M^m &
&\longrightarrow& & N^n\\
 & & \searrow {p_M} & &  \swarrow{p_N} & \\
& & & B^b & & 
\endmatrix \leqno (1.1)$$


\noindent Here $B^b,M^m$ and $N^n$ are  smooth connected manifolds of the indicated dimensions, without boundary, $B^b$ and $M^m$ being compact. Moreover $p_M$ and $p_N$ are smooth fibre  maps with fibres $F_M$ and $F_N$, resp. The (continuous) maps $f_1,f_2,f,...$ as well as homotopies between them are always assumed to be fibre-preserving (so that e.g. $p_N\circ f=p_M$ ); we also call them  {\it maps and homotopies over $B$} and write $f\sim_B f'$ \ 
if $f,f'$ are homotopic in this sense. From now on we will drop the superscript which denotes the dimension of the manifold, unless this simplification is going to cause some confusion.

\bigskip

\noindent{\bf Question 1.2.} \ Can the coincidence locus 

$$C(f_1,f_2):=\{x\in M \ | \  f_1(x)=f_2(x)\} $$

\bigskip
\noindent be made empty by suitable homotopies of $f_1$ and $f_2$ over $B$?  (If $f_1$ and $f_2$ can be deformed away from one another in this way we say that the pair $(f_1,f_2)$ is   {\it loose over $B$} or, shortly, {\it $B$-loose}).

More generally, we would like to estimate the {\it minimum number of pathcomponents}

$$MCC_B(f_1,f_2):=min\{\#\pi_0(C(f_1',f_2')) \  |  \  f_i\sim_B f_i', \ i=1,2\} \leqno (1.3)$$

\bigskip

\noindent of coincidence subspaces in $M$, achieved by suitable deformations of $f_1$ and $f_2$ over $B$.

For this purpose we study the geometry of the map


$$\matrix
(f_1,f_2)\ \ \colon  \ \  M\longrightarrow& N\times_B N &:=\{(y_1,y_2)\in N\times N \ | \ p_N(y_1)=p_N(y_2)\}\\
&\bigcup&\\
&\Delta &:=\{(y_1,y_2)\in N\times N \ | \ y_1=y_2\} \ . \ \ \ \ \ \ \ \ \ \ \    
\endmatrix \leqno (1.4)$$

\bigskip

After small deformations of $f_1$ and $f_2$ over $B$ this map is smooth and transverse to the diagonal $\Delta$ so that the coincidence locus

$$C=C(f_1,f_2)=(f_1,f_2)^{-1}(\Delta)=\{x\in M \ | \  f_1(x)=f_2(x)\} \leqno (1.5)$$

\bigskip

\noindent is an $(m-n+b)-$dimensional smooth submanifold of $M$.

Moreover the tangent map of $(f_1,f_2)$ induces an isomorphism of the normal bundles

$$\bar g_B^{\#}:\nu(C,M)\cong((f_1,f_2)|C)^{*}(\nu(\Delta, N\times_BN))\cong f_1^*(TF(p_N))|C \ ; \leqno (1.6)$$

\bigskip

\noindent here $TF(p_N)$ denotes the tangent bundle along the fibres of $p_N$.

A third important coincidence datum is the lifting


$$\matrix
& & E_B(f_1,f_2)&:=\{(x,\theta)\in M\times P(N) \ | \ p_N\circ\theta\equiv p_M(x),\\
& {\tilde{g}_B}\nearrow & \downarrow{pr}&\theta(0)=f_1(x),\ \theta(1)=f_2(x)\}\\
C& \longrightarrow& M \\
&{g=incl}& &
\endmatrix \leqno (1.7)
$$


\bigskip

 \noindent defined by $\tilde g_B(x)=(x,constant \  path \ at \ f_1(x)=f_2(x)$). Here $P(N)$ (and pr resp.) 
denote the space of all continuous paths    $\theta:[0,1] \to N,$ with the compact-open topology (and the obvious projection, resp.; compare \cite{10}, (5)).

The bordism class

$$\omega^{\#}_B(f_1,f_2)=[(C\subset M, \tilde g_B, \bar g_B^{\#})]\leqno (1.8)$$

\bigskip

\noindent of the resulting triple of the coincidence data (1.5)-(1.7) (which keeps track of the embedding of $C$ in $M$ 
and of its nonstabilized normal bundle) is independent of our choice of small deformations. It is our strongest  
obstruction to making the pair $(f_1,f_2)$ loose over $B$.
In certain settings (e.g. if $N=B \times S^ {n-b})$ it yields a complete homotopy classification for maps over $B$.
 However, this  (``strong") $\omega$-invariant is often hard to compute.

   The {\it stabilized} version

$$\tilde \omega_B(f_1,f_2)=[(C,\tilde g_B, \bar g_B)]\in \Omega_{m-n+b}(E_B(f_1,f_2);\tilde \varphi)\leqno (1.9)$$

\bigskip

\noindent is much more manageable. It forgets about the map $g:=pr\circ \tilde g_B$ (cf. 1.7) being an embedding, retains only the {\it stable} vector bundle isomorphism

$$\bar g_B: TC\oplus g^*(f_1^*(TF(p_N)))\oplus   \Bbb  R^k \cong    g^*(TM) \oplus \Bbb R^k, \ \ \ \ k>>0, \leqno (1.10)$$

\bigskip

\noindent (compare 1.6) and lies in the normal bordism group of singular $(m-n+b)-$manifolds in $E_B(f_1,f_2)$, with coefficient bundle

$$\tilde \varphi:= pr^*(f_1^*(TF(p_N))-TM)=pr^*(f_1^*(TN)-p_M^*(TB)-TM) \leqno (1.11)$$

\bigskip

\noindent(compare e.g. \cite{7},   2.1). The path space $E_B(f_1,f_2)$ and the resulting normal bordism group 
depend on the maps $f_1, f_2$, but  homotopies induce group isomorphisms  which preserve the $\tilde \omega_B$-invariants
(compare \cite{9}, 3.3). Therefore 
$\tilde \omega_B(f_1,f_2)$  vanishes    if $f_1$ and $f_2$ can be deformed to become coincidence free. In a suitable "stable dimension range'' the converse holds.

\proclaim{Theorem 1.1} Assume that 
$$m<2(n-b)-2.$$
\noindent Then a   pair $(f_1,f_2)$ is  loose over $B$ if and only if $\tilde \omega_B(f_1, f_2)=0.$
\endproclaim



 In the  proof  \ (outlined in the section 2 below) the path-space $E_B(f_1,f_2)$ plays a significant r\^ole: the lifting $\tilde g_B$ (cf. 1.7) allows us to construct 
the homotopies which deform $f_1$, $f_2$ away from one another. Quite generally $E_B(f_1,f_2)$ is a very interesting space with a rich topology. Already its decomposition into pathcomponents leads to the fibre theoretical  analogue of the (algebraic) Reidemeister equivalence relation (on $\pi_1(F_N))$  and to the corresponding notion of the {\it Nielsen numbers}

$$  N_B(f_1,f_2) \leq N_B^{\#}(f_1,f_2) \leq MCC_B(f_1,f_2).\leqno (1.12)$$

\bigskip

These are nonnegative integers counting  the path-components of $E_B(f_1,f_2)$ which contribute non-trivially to $\tilde \omega_B(f_1,f_2)$ and $\omega^{\#}_B(f_1,f_2)$,
 resp. (for details see section 4 below). Clearly these Nielsen numbers form lower bounds for the minimum number $MCC_B(f_1,f_2)$ (cf. 1.3); in particular, they are simple numerical looseness obstructions. Moreover the Nielsen numbers are  obviously smaller or equal to the {\it geometric Reidemeister  number}

$$\#R_B(f_1,f_2):=\#\pi_0(E_B(f_1,f_2)) \leqno (1.13)$$

\bigskip

\noindent (i.e. the number of path-components of the space $E_B(f_1,f_2)$, cf. 1.7; its relation to the 
classical  (algebraic) Reidemeister number will be explained in section 3). 
   
Another simplification of our $\tilde \omega_B$-invariant forgets about the path-space \break
$E_B(f_1,f_2)$ and the lifting $\tilde g_B$ altogether and keeps track only of the inclusion $g:C \subset M$ (as a continuous map) and of the description  of the stable normal bundle of $C$ given by (1.10). We obtain the normal bordism class

$$\omega_B(f_1,f_2)=[(C,g,\bar  g_B)]\in \Omega_{m-n+b}(M;\varphi) \leqno(1.14)$$

\bigskip

\noindent where 

$$\varphi:=f_1^*(TF(p_N))-TM=f_1^*(TN)-p_M^*(TB)-TM \leqno (1.15)$$

\bigskip

\noindent (compare 1.11). Homotopies of $f_2$ yield bordant triples of coincidence data $(C,g, \bar g_B)$ and hence the same
 $\omega_B-$invariants.

\bigskip

 {\bf Special case 1.16 (trivial base space):} If the base space $B$ consists of a single point we drop the subscript $B$ from our notations and obtain the invariants \ $\omega^{\#},\  N^{\#}, \ \tilde \omega, \  N$ and $\omega$ discussed  in 
\cite{8}, \cite{9} and    \cite{10}.  (For further literature concerning this special case see e.g. \cite{2}, \cite{3},
\cite{4},  \cite{5}, \cite{11} and    \cite{12} as well as the references listed there).

\bigskip

 {\bf Special case 1.17 (trivial target fibration):} If the target fibration is a product, $N=B\times F_N$, we may write $f_i=:(p_M, f_i'), \ i=1,2$.  Then the $\omega_B^{\#}-, \ \tilde \omega_B-$ and $\omega_B-$invariants of $(f_1,f_2)$ are related to the corresponding (unfibered) invariants of $( f_1', f_2')$ via bijections (which preserve 0); in particular 

$$N_B^{\#}(f_1,f_2)=N^{\#}( f_1', f_2') \ \   \hbox{and} \ \    N_B(f_1,f_2)= N( f_1', f_2').$$

\bigskip

 {\bf Special case 1.18 (fixed points):}  If the two fibrations coincide and $f_1$ is the identity map $id$, then $C(id,f)$ is the fixed point locus of $f$, the coefficient bundles $\tilde \varphi$ and $\varphi$ are the pullbacks of the virtual vector bundle $-TB$ under  projections, and our $\omega-$invariant can be weakened further to yield the bordism class

$$p_{M*}(\omega_B(id,f))=[(C @> {p_M|}>> B, \bar g_B: TC @ >{\cong }> {stably} > (p_M|)^*(TB))] \ \in \ \Omega_b(B; -TB)\leqno (1.19)$$

\newpage


\noindent This procedure neglects the "vertical" aspects of our fixed point data.

 On the other hand A. Dold \cite{1} has defined a fixed point index $I^h(f)$ of $f$ for every multiplicative  generalised cohomology theory $h$ with unit. In view of the universality property of stable cohomotopy theory the  strongest ("universal")  version of Dold's index
 takes the form 

$$I(f)\in \pi_{stable}^0(B^+)=\varinjlim [\Sigma^k B^+, S^k]; \leqno(1.20)$$

\bigskip

\noindent and actually classifies certain "horizontal" fixed point phenomena (cf. \cite{1}, theorem 4.3); here $B^+$ denotes  the space $B$ with a disjointly  added point.

Note that the Pontrjagin-Thom  procedure yields a canonical isomorphism

$$PT \ \ : \ \  \pi^0_{stable}(B^+) \ \ @> \cong >> \ \  \Omega_b(B;-TB) \leqno (1.21)$$

\bigskip
\noindent (which will be described in section 5 below).



\proclaim{Theorem 1.2}  For every map $f:M \to M$ over $B$ 

$$I(f)=(PT)^{-1}(p_{M*}(\omega_B(id,f))).$$  
\endproclaim        

\bigskip
 
\noindent  The proof will be given in section 5 below.
\qed \quad

 \bigskip

As in illustration of our notions and methods we calculate the minimum number $MCC_B$ (as well as the Reidemeister and the Nielsen   numbers) and the $\omega_B$-invariant for      all pairs of $B-$maps involving the torus and$\slash$or the Klein bottle over $B=S^1$. Note that this is way  outside of the stable dimension range discussed in  theorem 1.1.

\medskip

{\bf Example 1.22 ($S^1-$bundles over $S^1$):} Let $M,N$ be (possibly different) fibre spaces over $S^1$ with fibre $S^1$. 
Thus $M$ (and also $N$) is either the torus
$$T=S^1\times S^1=I\times S^1/(0,z)\sim(1,z),z\in S^1 \leqno (1.23)$$ 

\bigskip

\noindent or the Klein bottle

$$K=I\times S^1/(0,z)\sim(1,\bar z), z\in S^1 \leqno (1.24)$$

\bigskip
\noindent with the standard  projection to \ \   $I/0\sim 1=S^1$. We define two sections $s_{\epsilon}$, $\epsilon=\pm 1$, by 

$$s_{\epsilon}([t]) \ = \ [(t,\epsilon)].\leqno (1.25)$$

\bigskip

 Given a map $f:M \to N$ over $S^1$ we have two well defined numerical invariants:

$$   q(f):=(degree \ \ of \ \  f|:F_M\to F_N) \ \in \  \Bbb Z \ \ \ \ \leqno(1.26)$$

\noindent (this vanishes if $M\ne N$); \ \ \ \ \  and

$$    r(f):= degree \ \ of \ (B=S^1 @>  {f\circ s_{+1}}>>  N=[0,1]\times S^1/\sim \ \ \ \longrightarrow S^1);\leqno (1.27)$$

\noindent this lies in $\Bbb Z$ (and in $\Bbb Z_2$, resp.) if $N=T$ (and if $N=K$ {\it and} $f$ preserves the base point $[(0,1)]$, resp.); this number measures roughly how   often the section $f\circ s_{+1}$ (assumed to be base point preserving if $N=K$) ``winds around the fiber in $N$''. A base point free description  of $r(f)$ in the case $N=K$ is as follows: $r(f)$ equals the mod 2 integer  
0 (and 1, resp.) if   $f\circ s_{+1}$ is homotopic (through sections in $K$) to $s_{+1}$ (and to $s_{-1}$, resp.).

\bigskip

Returning to the base   point free  setting we obtain:


\proclaim{Proposition 1.1} Two maps $f, \hat f: M \to N$ over $S^1$ are  homotopic  over $S^1$ if and only if 
$q(f)=q(\hat f)$ and $r(f)=r(\hat f)$.
\endproclaim

Thus each homotopy class can be represented by a map in a rather natural standard form (enjoying constant angular velocities both along each fibre and for $f\circ s_{+1})$. This is very helpful when we analyse coincidence data.

 Now consider any two maps $f_1,f_2: M \to N$ over $S^1$ and put 
$$q:=q(f_1)-q(f_2) \ \  \ \ \hbox{and} \  \ \ \ r:=r(f_1)-r(f_2)  \leqno(1.28)$$

\noindent (compare 1.26 and 1.27).

\proclaim{Theorem 1.3} The minimum number $MCC_B(f_1,f_2)$ is equal to the Nielsen numbers 
$  N_B(f_1,f_2)$ and $N^{\#}_B(f_1,f_2)$ $($and   also to $\#R_B(f_1,f_2)$ whenever this Reidemeister number is 
finite  $)$.
More precisely:\\
(i) Assume $N=S^1\times S^1$. Then $(q,r)\in \Bbb Z\times \Bbb Z$ and we have: 
$$gcd(q,r)=MCC_B(f_1,f_2)=N_B(f_1,f_2)=\#R_B(f_1,f_2) \ \  if \ \  (q,r)\ne (0,0);$$ 
\noindent  $$0=MCC_B(f_1,f_2)=N_B(f_1,f_2)\ne \# R_B(f_1,f_2)=\infty \ \ \ \   if \ \  (q,r)=(0,0).$$

 In particular, the pair $(f_1,f_2)$
is loose over $B$ if and only if $f_1 \sim_B f_2$.\\


\newpage

(ii) Assume $N=K$. Then  $(q,r)\in \Bbb Z\times \Bbb Z_2$ and we have: if  $q\ne 0$:

$$ MCC_B(f_1,f_2)=N_B(f_1,f_2)=\#R_B(f_1,f_2)= \cases 
|q|/2 &\ \hbox{if}\ q \ even, \ r=1\\
[ \ |q|/2 \ ]+1 &\ \hbox{else};
\endcases
$$

$$ if \ q=0: MCC_B(f_1,f_2)=N_B(f_1,f_2)=\cases 
0 & \hbox{if}\ r\ne 0 \\
1 & \hbox{if} \ r=0
\endcases  \ \ \ 
\ne \#R_B(f_1,f_2)=\infty.
$$

 In particular,   the pair $(f_1,f_2)$
is loose over $B$ if and only if it consists of two ``antipodal`` maps, i.e. \ $-f_1\ \ \sim_B \ \ f_2$.   
\endproclaim

Note that here the value of the Nielsen number is always 0 or 1 or the Reidemeister number. A similar result in an entirely different setting was proved in \cite{12}, theorem 1.31.

Clearly, in all of example 1.22 the $\tilde \omega_B$-invariant is a complete looseness obstruction. Actually, already its
    weaker version $\omega_B(f_1,f_2)\in \Omega_1(M; \varphi)$(cf. 1.14) allows us to distinguish maps up to homotopy over $S^1$.

\proclaim{Theorem 1.4} Let  $(M,N)$ be any of the four combinations of $S^1-$bundles over $S^1$ and let 
$f_1,f_2:M \to N$ be maps over $S^1$. Then there are  canonical isomorphisms which describe $\Omega_1(M; \varphi)$ $($and correspondingly $\omega_B(f_1,f_2)$$)$ as $($an element of$)$ a direct sum of three groups, as 
follows   $($compare proposition 1.1 and theorem  1.3$)$\\





\settabs5\columns
\+$(M,N)$&$\Omega_1(M;\varphi)$& &\multispan3$\omega_B(f_1,f_2)$&\omit\cr
\+$(T,T)$&$\Bbb Z\oplus \Bbb Z\oplus \Bbb Z_2$&$q$&$r$&\omit\cr 
\+$(K,K)$&$\Bbb Z\oplus \Bbb Z_2\oplus \Bbb Z_2$&$q$&$r+1+\rho_2(q)$&$\rho_2(N_B(f_1,f_2))$\cr 
\+$(K,T)$&$0\oplus \Bbb  Z\oplus \Bbb Z_2$&$q$&$r$&(compare \cr 
\+$(T,K)$&$0\oplus \Bbb  Z_2\oplus \Bbb  Z_2$&$q$&$r+1$&theorem 1.3)\cr

\bigskip

 Here  $\rho_2:\Bbb Z \to \Bbb Z_2$ denotes reduction mod 2.

In particular, for every map $f: M \to N$ over $S^1$ the ``fibred degree" $($or ``root invariant"$)$ $\omega_B(f,s_{+1}\circ p_M)$ determines $q(f)$ and $r(f)$ and hence the homotopy class $($over $S^1$$)$ of $f$.
\endproclaim

{\it Remark 1.29.} In view of     proposition 1.1 the homotopy class of $f$ is already determined by the first two components of $\omega_B(f, s_{+1}\circ p_M)$ or, equivalently, by 

$$\mu(\omega_B(f,s_{+}\circ p_M))\in H_1(M;\tilde \Bbb Z_{\varphi})$$

\newpage

\noindent where $\mu$ denotes the Hurewicz homomorphism into the first homology group of $M$ with integer coefficients (which are twisted like the orientation line bundle of $\varphi$). This is a very special phenomenon, related to the fact that both the torus and the Klein bottle are $K(\pi, 1)'s$. For general $M$ and $N$ the methods of singular homology theory are often far too weak, and the full power of our approach   (based on normal bordism theory and the   pathspace 
$E_B(f_1,f_2)$) yields better results.

{\it Remark 1.30.} Consider a selfmap of the torus $T$ or the Klein bottle $K$ over $S^1$. Dold's fixed point  index 
 \cite{1} in  its strongest form lies in 

$$\pi^0_{stable}((S^1)^+) \ \   \cong \ \    \Omega_1^{fr}(S^1) \ \   \cong \ \ \Bbb Z\oplus \Bbb Z_2$$

\bigskip

\noindent and captures precisely the first and third components    of $\omega_B(id, f)$ or, 
equivalently,   $\pm q=\pm(deg(f|F)-1)$ as well as the Nielsen number $N_B(f,id)$, taken mod  2. However, it looses all information about the characterictic winding number $r$ which $-$ together with $q$ $-$    determines $f$ and which measures the ``vertical" aspect of the generic fixed point circles.

\head 2. The $\omega$-invariants \endhead

  Given maps $f_1,f_2: M \to N$ over $B$, the definition of $\omega^{\#}_B(f_1,f_2)$, 
 $\tilde \omega_B(f_1,f_2)$ and  $\omega_B(f_1,f_2)$  (as outlined in the introduction) is completely analoguous to the definition  (given in \cite{10} and \cite{9}) of the corresponding invariants for ordinary maps between manifolds (or, equivalently, for maps over $B=\{point\}$). Therefore many of the notions, methods and results of the ordinary (fibration free) coincidence theory allow a straightforward generalization to the setting of  fibre preserving maps.  

 In particular, the  proof of theorem 1.1 proceeds in direct analogy to the proof of theorem 1.10 in (\cite{9}, pp 213 and 223-224): we just  have to replace $N\times N$ by $N\times_BN$. Our 
(``stable'') dimension condition means that the dimension of $C(f_1,f_2)$,  augmented by 2, is strictly smaller than the codimension in $M$. Hence here \     $\tilde \omega_B(f_1,f_2)$ is precisely as strong as    $\omega^{\#}_B(f_1,f_2)$, and
 $N_B(f_1,f_2)=N_B^{\#}(f_1,f_2)$; nulbordism data can be realized by a suitably embedded manifold in $M \times I$ with a nonstabilized description of its normal bundle and, above all, without new coincidences occurring in its shadow  (cf. 
\cite{9}, p. 224).

\bigskip

{\it Remark 2.1.} The interested reader may check when the methods of \cite{9}, 1.10 and 4.7, can be generalized to yields the full equality $MCC_B(f_1,f_2)=N_B(f_1,f_2)$ (``Wecken theorem'').

\bigskip

 Next let us consider the special case where the target fibration is trivial. Given maps over $B$, 
$$f_i=(p_M, f'_i): M \to N=B\times F_N  , \ i=1,2;$$ 
\noindent we see that $E_B(f_1,f_2)$ can be identified with the path-space $E(f_1',f_2')$ discussed in  \cite{9} and  \cite{10}. Thus the $\omega$-invariants and Nielsen numbers {\it over $B$} of the pair $(f_1,f_2)$ are equal to the corresponding {\it ordinary}   (unfibred) invariants of $(f_1', f_2')$.

\head 3. The algebraic  Reidemeister classes    over $B$ and the space $E_B(f_1,f_2)$ \endhead

 In this section we fix maps $f_1,f_2: M \to N$ over $B$. We  will give an algebraic description of the  (geometric) Reidemeister set 
$\pi_0(E_B(f_1,f_2))$ (compare 1.13). This generalizes and refines the classical approach. As an application we will compute  Reidemeister numbers for maps into the Klein bottle.

 Choose a coincidence point $x_0 \in C(f_1,f_2)$ (if it  does not exist, the pair $(f_1,f_2)$ is loose and our initial question 1.2 needs no further answer). Put \break $y_0:=f_1(x_0)=f_2(x_0)$ and let $F_N \subset N$ be the fibre over $b_0:=p_M(x_0)$.

 Using  homotopy lifting extension properties (compare \cite{13}, I.7.16) of the (Serre) fibration  $p_N$ we construct a well defined operation

$$*_B \ : \  \pi_1(M, x_0)\times \pi_1(F_N, y_0) \longrightarrow \pi_1(F_N, y_0)  \leqno (3.1)$$
\bigskip
\noindent as follows. Given loops $c:(I, \partial I) \to (M,x_0)$ and $\theta : (I, \partial I) \to (F_N, y_0)$, lift the
 homotopy  

$$ h: I \times I \to B, \ \ \ \ \ \  h(s,t) \ : \ = p_M \circ c(s), \leqno (3.2)$$
\bigskip 
\noindent to a map $\tilde h: I \times I \to N$ such that 
$$\tilde h (0,t)=\theta(t) \ \ ,\  \  \tilde h(s,0)=f_1\circ c(s)\ \ ,\ \  \tilde h(s,1)=f_2\circ(s) \leqno (3.3)$$
\bigskip
\noindent for all $s,t \in I$. Then the loop $\theta'$ defined by $\theta'(t):=\tilde h(1,t)$ lies entirely in $F_N$. Due to the very special form of $h$ (cf. 3.2) the homotopy class $[\theta']$ of $\theta'$ in $F_N$ (and {\it not} just in $N$) depends only on the homotopy classes of $c$ and $\theta$. We put 

$$ [c]*_B[\theta] \ : \ = \ [\theta']. \leqno (3.4)$$
\bigskip

\proclaim{Definition 3.1} Two elements $[\theta]$, $[\theta']$ $\in \pi_1(F_N,y_0)$ are called  \ {\rm Reidemeister equivalent over $B$}   if there exists $[c]\in \pi_1(M,x_0)$ such that 
$[c]*_B[\theta]=[\theta']$.  

The {\rm algebraic Reidemeister set $R_B(f_1,f_2, x_0)$} is the resulting  set of  equivalence classes (i.e. of orbits of the group action  $*_B$ of
 $\pi_1(M,x_0)$ on (the set)  $\pi_1(F_N, y_0)$). 

Its cardinality is called {\rm Reidemeister number of $f_1,f_2$ over $B$}.
\endproclaim

There is also the classical group action (without any reference to $B$)

$$ * \ : \ \pi_1(M, x_0) \times \pi_1(N, y_0) \longrightarrow \pi_1(N, y_0) \leqno (3.5)$$
\bigskip

\noindent determined by the induced homomorphisms $f_{j*}: \pi_1(M, x_0) \to \pi_1(N, y_0),$ $j=1,2$, \ i.e.



$$[c] \ * \ [\theta]\  \ : \ = \ \ f_{1*}([c])^{-1} \ \cdot \ [\theta] \ \cdot  \  f_{2*}([c]) \leqno(3.6) $$

\bigskip
\noindent for $[c] \in\pi_1(M,x_0)$ and $[\theta]\in\pi_1(N,y_0)$ (compare e.g. \cite{9}, 2.1).

In view of the boundary conditions (3.3) of the lifting $\tilde h$ (which takes its value in $N$  and, in general , not already in $F_N$) we see that 

$$[c]* i_*([\theta]) \ = \ i_*([c]*_B[\theta]) \leqno (3.7)$$

\bigskip

\noindent for all $[c]\in \pi_1(M,x_0), \ [\theta]\in\pi_1(F_N,y_0);$ here  $i:F_N \to N$ denotes the inclusion. In particular, the standard action $*$ (cf. 3.6) restricts to an action of $\pi_1(M,x_0)$ on $i_{*}(\pi_1(F_N,y_0)).$ In general this yields a coarser equivalence relation than the one defined by  our action \ $*_B$ \ (e.g. when 
\  \ $p_M=p_N : S^{2k+1} \to \Bbb CP(k)$ \ , \ $k\geq 1$,  is the Hopf fibration, then $R_B(f_1,f_2,x_0)=\pi_1(S^1)\cong \Bbb Z$, but $i_*(\pi_1(F_N))=0$). However, if $i_*$ is injective (e.g. when $\pi_2(B)=0$) then  (3.7) can be used to compute

$$R_B(f_1,f_2,x_0)=\pi_1(F_N,y_0)/\sim *_B \approx \ i_*(\pi_1(F_N,y_0))/\sim*. \leqno(3.8)$$

\bigskip
\noindent In particular, when $B=\{b_0\}$ and hence $F_N=N$,  our definition of an algebraic Reidemeister set coincides with the usual one.

More general injectivity criteria for $i_*$ may be extracted from the exact sequence

$$  \cdots \to \pi_2(B, b_0) \to \pi_1(F_N,y_0) \ @>{i_*}>>\pi_1(N,y_0) @> {p_{N*}}>>\pi_1(B, b_0). \leqno (3.9)$$

\bigskip

 Next let us compare our algebraic and geometric Reidemeister sets (cf. definition 3.1 and (1.13)). By definition $E_B(f_1,f_2)$ is the space of pairs $(x,\theta)$ where $x$ is a point  in  $M$ and  $\theta$ is a path in $N$ from $f_1(x)$ to $f_2(x)$ which stays entirely in one fibre of $p_N$. In  view of the very special form of the homotopy $h$ (cf. 3.2) its lifting $\tilde h$ determinies a path
$$s\in I \ \longrightarrow \ ((c(s), \tilde h(s,-)) \ \in \ E_B(f_1,f_2) $$
\bigskip
\noindent joining $(x_0, \theta)$ to $(x_0,\theta')$. Actually every other path in $E_B(f_1,f_2)$ which starts and ends in the fibre $pr^{-1}(\{x_0\})=\{x_0\}\times \Omega(F_N,y_0)$ of $pr$ \  (cf. 1.7) can be obtained in this way from some  lifted
homotopy $\tilde h$ as in (3.2), (3.3). In other words, two classes $[\theta], \  [\theta']\in  \pi_1(F_N, y_0)$ are Reidemeister equivalent over $B$  if and only if $(x_0,\theta)$ and $(x_0,  \theta')$ lie in the same path-component
 of $E_B(f_1,f_2)$.

 Thus the map 
$$R_B(f_1,f_2,x_0) \longrightarrow \pi_0(E_B(f_1,f_2)) \ ,  \leqno (3.10)$$
\bigskip
\noindent which is induced by the fibre inclusion  \ 
  $\Omega(F_N,y_0)\approx pr^{-1}(\{x_0\}) \  \subset \ E_B(f_1,f_2)$   \\
and by the resulting map   
$$\pi_1(F_N, y_0)=\pi_0(\Omega ( F_N, y_0)) \longrightarrow \pi_0(E_B(f_1,f_2)) \ , $$ 
\bigskip

\noindent is injective. It is also onto. Indeed, given any point $(x, \theta)$ of $E_B(f_1,f_2)$, we can pick a path in $M$ from $x$ to $x_0$ and lift it to a path in $E_B(f_1,f_2)$ which joins $(x, \theta)$ to some point in $pr^{-1}(\{x_0\})$.

We have showed

\proclaim{Theorem 3.1}
For every pair $f_1,f_2: M \to N$ of maps over $B$ and for every choice $x_0\in C(f_1,f_2)$ and $y_0=f_1(x_0)=f_2(x_0)$ of base points there is a canonical bijection 

$$R_B(f_1,f_2,x_0) \ \approx \ \pi_0(E_B(f_1,f_2))$$
\bigskip
\noindent between the algebraic and geometric Reidemeister sets.
 \endproclaim

\proclaim{Corollary 3.1} The Reidemeister number depends only on the (base point free) homotopy classes of $f_1$ and $f_2$ over $B$.
\endproclaim
Indeed, any pair of homotopies $f_1 \sim f_1'$, $f_2 \sim f_2'$ over $B$ induces a fibre homotopy equivalence $E_B(f_1,f_2) \sim E_B(f_1', f_2')$ over $M$ (compare \cite{9}, 3.2). \quad \qed

\bigskip

{\it EXAMPLE 3.11} (Maps into the Klein bottle).  We illustrate the previous discussion by a calculation which we will need in the  proof of theorem 1.3.

\proclaim{Proposition 3.1} Consider maps $f_1,f_2: M \to K$ over $S^1$ where $M$ is the torus $T$ or the Klein bottle $K$ (and use the notations (1.26)-(1.28)).
 
If $M=T$ or $q=0$, then the Reidemeister number $\#R_B(f_1,f_2)$ is infinite.

 If $M=K$ and $q\ne 0$, then

$$\#R_B(f_1,f_2)= \cases 
|q|/2 &\ \hbox{if}\ \ q\equiv 0(2), \ r\ne 0;\\
[\ |q|/2\ ]+1 &\ \hbox{else}.
\endcases
$$

\endproclaim

\demo{Proof} In view of corollary 3.1 we may assume that $f_1$ and  $f_2$ map \break $x_0=[(0,1)]$ to  $y_0=[(0,1)]$ (cf. (1.23) and (1.24)).  Let us use these base 
points for computing the algebraic Reidemeister set. Then $\pi_1(M)$ (and $\pi_1(K)$, resp.) is generated by 

$$a_M  :=  i_{M*}(g) \ \ \ \ \hbox{and} \ \ \ \  b_M  :=  s_{+1*}(g) $$

\bigskip

\noindent (and \ by \  $a:=  i_*(g)$   and    $b  :=  s_{+1*}(g)$,  resp.)  where \ $i_M, \  i, \  s_{+1}$ \ denote fibre inclusions and the section defined in (1.25); $g$ is the standard generator of $\pi_1(S^1)$.

 Since $\pi_2(S^1)$ vanishes, $i_*$ is injective and we have to evaluate only the standard action (3.6) of $\pi_1(M)$ on $\pi_1(F_N)\cong \Bbb Z.$ Given $k\in  \Bbb Z$, we obtain 

\bigskip

\noindent (3.12) \ \ \ \ \ \ \ \ \ \ \ \  $a_M  *a^k \  = \ a^{k-(q(f_1)-q(f_2))} \ = \ a^{k-q} \  ; $

\bigskip

$ b_M *a^k \  = \ f_{1*}(b_M)^{-1}\cdot a^k\cdot f_{2*}(b_M)\ = \ b^{-1} \cdot a^{-r(f_1)} \cdot a^k \cdot a^{r(f_2)}\cdot b \ = \ a^{r(f_1)-r(f_2)-k}$

\bigskip

\noindent  where we consider $r(f_j)\in \{0,1\}$ as an {\it integer} so that $f_{j*}(b_M)=a^{r(f_j)}\cdot b \ , \ j=1,2$ (compare (6.1)). Therefore we can 
interpret \  $\pi_1(F_N)/\sim *_B $ \ \ (cf. (3.8)) as the orbit set of the involution $\iota$ on $\Bbb Z/ q \Bbb Z$ defined by 

$$\iota([k]) \ := \ [r(f_1)-r(f_2)-k],\ \  \ \ \ [k]\in \ \Bbb Z/ q \Bbb Z. \leqno (3.13) $$

\bigskip

In particular,  its cardinality is infinite if $q=0$. This is e.g. always the case when $M=T$, since the map $f_j| \ : \ F_M \to F_N=S^1$ is freely homotopic  to its own complex conjugate and hence has degree $q(f_j)=0$, $j=1,2$.

 For the remainder of the proof it suffices to consider the case where $M=K$ and $q>0$. Then 
$$\#R_B(f_1,f_2,x_0) \  = \  (q+\# Fix(\iota))/2. \leqno (3.14)$$ 

\bigskip
 Clearly the fixed point set $Fix(\iota)$ of $\iota$ consists just of the solutions of the linear equation
$$2[k] \ = \  [r(f_1)-r(f_2)]$$

\bigskip
\noindent in  $\Bbb Z/q\Bbb Z.$ Therefore it is easy to see that

$$\#Fix(\iota)= \cases 
1 &\ \hbox{if}\ q\ is \ odd;\\
2 &\ \hbox{if} \ q\equiv 0(2), \  r(f_1)=r(f_2);\\
0 &\ \hbox{if} \ q\equiv 0(2), \  r(f_1)\ne r(f_2).
\endcases
$$

\bigskip
\noindent In view of (3.14) this completes the proof.
\quad \qed
\enddemo

\head 4.  Nielsen  coincidence  classes over $B$ \endhead

 In this section we extend J. Jezierski's notion of Nielsen classes over $B$ (cf. \cite{6}) in the obvious way from fixed points to coincidences of maps $f_1$, $f_2$ over $B$. The resulting decomposition of the coincidence set turns out to correspond precisely to the decomposition of the space 
$$ E_B(f_1,f_2) \ = \ \bigcup _ {A \in \pi_0(E_B(f_1,f_2))}  A \leqno (4.1)$$

\bigskip  
\noindent into path-components and yields the description of 

$$\tilde \omega_B(f_1,f_2) \ = \ \{(\tilde \omega_B(f_1,f_2))_A\} \ \in \ \Omega_*(E_B(f_1,f_2); \tilde \varphi) \ = \ \oplus_{A} \ \Omega_*(A; \tilde \varphi|A) \leqno (4.2) $$

\bigskip 
\noindent as a direct sum. We will discuss the Nielsen number

$$N_B(f_1,f_2)\ \ :=\  \ \#\{A\in \pi_0(E_B(f_1,f_2)) \ | \ (\tilde \omega_B(f_1,f_2))_A \ne 0\} \leqno (4.3) $$

\bigskip
\noindent (which counts the nontrivial direct summands of $\tilde \omega_B(f_1,f_2))$ and its nonstabilized analogue

$$N_B^{\#}(f_1.f_2) \ := \ \#\{A \in \pi_0(E_B(f_1,f_2)) \ | \ (\omega_B^{\#}(f_1,f_2))_A \ne 0\}.
 \leqno (4.4) $$
\bigskip

\noindent In classical fixed point theory (where  $B$ consists of a single point) both definitions (4.3) and (4.4) just yield the  familiar notion of the Nielsen fixed point number.

\proclaim {Definition 4.1} Let $f_1,f_2 : M \to N$ be maps over $B$. Two coincidence points $x,x'\in C(f_1,f_2)$ are called {\rm Nielsen equivalent over $B$} if there exist
 a path $c: I \to M$ joining $x$ to $x'$, as well as a homotopy $\tilde h: I\times I \to N$ from $f_1\circ c$  to $f_2\circ c$ which keeps the end points fixed and such that for each $s\in I$ the whole image $\tilde h(\{s\}\times I)$ lies in the fibre of $p_N$ over $p_M\circ c(s).$
\endproclaim

\proclaim{Proposition 4.1} The coincidence points $x$ and $x'$ are Nielsen equivalent over $B$ if and only if the map
 $$ \tilde g_B \ : \ C(f_1,f_2) \longrightarrow E_B(f_1, f_2)$$

\noindent (defined by \ $\tilde g_B(x)=(x, constant \ path \  at \ f_1(x)=f_2(x))$ takes them into the same path-component
$A$ of $E_B(f_1,f_2)$. Therefore the Nielsen classes of $(f_1,f_2)$ over $B$ are just those inverse images $\tilde g ^{-1}(A), \  A\in \pi_0(E_B(f_1,f_2))$, which   are nonempty.
\endproclaim

\demo{Proof} (Compare also the proof of theorem 3.1). The data $(c, \tilde h)$ in the definition 4.1  represent just another way of describing a path in $E_B(f_1,f_2)$ from 
 $\tilde g_B(x)$ to $\tilde g_B(x')$. Indeed,  for every $s\in I$ the pair $(c(s), \tilde h(s, -))$ lies in  $E_B(f_1,f_2)$, since $\tilde h(s, - )$ is a path
 joining $f_1(c(s))$ to $f_2(c(s))$ in the fibre $p_N^{-1}(p_M(c(s)))$.
\quad \qed
\enddemo

\proclaim {Corollary 4.1} Each Nielsen class is open and closed in   $C(f_1,f_2)$.
\endproclaim

Indeed, it is not hard to see that each path-component $A$ is open and closed in $E_B(f_1,f_2)$.
\quad \qed

\bigskip

 We want to consider only those Nielses classes which survive (in some sense) all possible $B$-homotopies of $f_1$, $f_2$. We try to detect them with the help of our
 $\omega-$invariants.

 After a suitable approximation of $f_1$, $f_2$ the coincidence set $C$ is a clossed manifold, and so is each Nielsen class $ C_A  :=  \tilde g_B^{-1}(A), \ A\in \pi_0(E_B(f_1,f_2)).$ We call it {\it strongly essential}, and {\it essential}, resp., if the corresponding triple \break 
$(C_A, \ \tilde g_B|C_A , \ \bar g^{(\#)}|C_A)$ of restricted coincidence data is not nullbordant (in the nonstabilized, and stabilized sense, resp.). Define $N^{\#}_B(f_1,f_2)$
 and $N_B(f_1,f_2)$ to be the resulting numbers of (strongly) essential Nielsen classes.

 \proclaim{Theorem 4.1} For all maps $f_1,f_2: M \to N$ over $B$ we have:\\

i) \ the  Nielsen numbers $N^{\#}_B(f_1,f_2)$ and $N_B(f_2,f_1)$ depend only on the homotopy classes of $f_1$, $f_2$ over $B$;\\

ii) \  $N^{\#}_B(f_1,f_2)=N^{\#}_B(f_2,f_1)$ \ \ \ \ and \ \ \ \   $N_B(f_1,f_2)=N_B(f_2,f_1)$;\\

iii) \  $0 \leq  N_B(f_1,f_2) \leq N^{\#}_B(f_1,f_2) \leq MCC_B(f_1,f_2) <\infty$ \ \ \ and 

 \ \ \ \ \ \ \ \ \ \ \ \ \ \ \ \ \ \ \ \ \ \ \ \ \ \ \    $N^{\#}_B(f_1,f_2) \leq \#R_B(f_1,f_2)$;\\

iv) \ in classical fixed point theory (over $B$=point) both versions of our Nielsen numbers coincide with the   classical notion of the Nielsen numbers.

\endproclaim

 The proof proceeds as in \cite{9}, 1.9, and \cite{10}, 1.2.  \quad \qed

\bigskip

Unlike the $\omega_B$-invariants which lie in (possibly very complicated) bordism sets (varying with $f_1$, $f_2$) our Nielsen numbers are simple numerical  looseness obstructions. To what extend are they less powerful?

\proclaim{Proposition 4.2} $N_B(f_1,f_2)=0$         if and only if $\tilde \omega_B(f_1,f_2)=0.$
\endproclaim
 This follows from  the direct sum decomposition 4.2.

 It is not  clear  whether the corresponding statement holds for $N^{\#}_B$ and $\omega_B^{\#}$. If $N_B^{\#}(f_1,f_2)=0$, then  the Nielsen classes $C_A$ allow individual embedded nullbordisms in $M\times I$. But these may not fit together disjointly to yield  an {\it embedded} nullbordism for all of $C(f_1,f_2)$ (which is needed to show 
that $\omega^{\#}(f_1,f_2)=0$).

\head 5. Relation to Dold's index \endhead

In this section we study the special case 1.18 where the two   fibrations $p_M$ and $p_N$ coincide, $f_1$ is the identity map $id$ and we are interested in the fixed point behaviour of a map $f_2=f$ over $B$.

 We will see that our weakened normal bordism invariant $p_{M*}(\omega(id,f))$ determines the strongest  version of Dold's fixed 
point index (which generalizes the Lefschetz index, cf. \cite{1}).

 First let us describe the Pontrjagin-Thom isomorphism PT (cf. 1.21) which relates these invariants. Given a real number $R>0$, let 
$\bar D^k(R)$ (and $D^k(R)$, resp.) denote the compact (and open, resp.) ball of radius $R$ in euclidian space $\Bbb R^k$ and identify the quotient space
$$\bar D^k(R)/\partial \bar D^k(R)=\Bbb R^k/(\Bbb R^k-D^k(R))$$

\bigskip
\noindent with  the sphere $S^k=\Bbb R^k\cup \{\infty\}$ in the standard fashion. Moreover define

$$E_R^k:=B\times \bar D^k(R)\ \  \subset \ \   E^k:=B\times \Bbb R^k.\leqno (5.1)$$

\bigskip

Then we can interpret the suspension
$$\Sigma^kB^+=B\times S^k/(B\times\{\infty\})=E^k/(E^k-\overset{\circ}\to E_R^k)=E_R^k/\partial E_R^k \leqno (5.2)$$

\bigskip
\noindent as a one point-compactification of $\overset{\circ}\to E_R^k=B\times D^k(R)$. Now, given a map

$$u:(\Sigma^kB^+, \infty) \to (S^k, \infty)\ \ \ , \ \ \ k >>0,$$
\bigskip
\noindent up to homotopy, we may assume that $u|\overset{\circ}\to E_R^k$ is  smooth with regular value $0\in \Bbb R^k\subset S^k.$ Thus its 
inverse image $u^{-1}(\{0\})$  is a  smooth submanifold of $B \times \Bbb R^k$ whose normal bundle is trivialized via the tangent map of
 $u$. The resulting normal bordism class $[(u^{-1}(\{0\}), first \ projection, \bar g_B) ]$ is the value of $[u]\in \pi_{stable}^0(B^+)$ under the Pontrjagin-Thom isomorphism PT (cf. 1.21; compare 1.10 and 1.20).

{\it Proof of Theorem 1.2.} In view of the homotopy  invariance of Dold's index $I(f)$ (cf. \cite{1}. 2.9) we may assume that the map 
$(id, f): M \to M\times_B M$ is smooth and transverse to the diagonal $\Delta$. Then the fixed point set 
$$C=C(id, f)=(id,f)^{-1}(\Delta) \leqno (5.3)$$
\bigskip
\noindent is a smooth submanifold of $M$ with the description

$$\bar g_B^{\#}: \  \nu_1:=\nu(C,M) \cong TF(p_M)|C \leqno (5.4)$$

\bigskip 
\noindent of its normal bundle as in (1.6).

 For large $k$ there exists a smooth embedding $M \subset B\times \Bbb R^k$ over $B$ whose normal bundle $\nu_2:=\nu(M, B\times \Bbb R^k)$ can be identified with a  ``vertical" subbundle of $T(B\times \Bbb R^k)$ which, together with $TF(p_M)$, spans the tangent bundle along the fibres of $B \times \Bbb R^k \to B$. Let $\bar V$ (and V, resp.) be a corresponding compact (and open, resp.) tubular  neighborhood of $M$ in $B \times \Bbb R^k$ and consider the composite map

$$\hat f \ \ :\ \  \bar V@>projection >>M @>f >> M\subset B\times \Bbb R^k \leqno (5.5)$$ 

\bigskip
\noindent over $B$. Clearly its fixed point set is also equal to $C$ (cf. 5.3). Hence there exists a radius $R>0$ such  that $v-\hat f(v)\notin B \times D^k(R)$ for every $v\in \partial \bar V.$ Moreover we can pick a radius $\rho>0$ such that the space $\bar V$ lies in $B\times \bar D^k(\rho).$ Collapsing its complement and using (5.1) and (5.2) we obtain the composite map

$$\hat u:\Sigma^kB^+=E_{\rho}^k/\partial E_{\rho}^k\to \bar V/\partial \bar V@>{id-\hat f}>>E^k/(E^k-\overset{\circ}\to{E}_ R^k)=\Sigma^kB^+. \leqno (5.6)$$ 
\bigskip

 Now, according to \cite{1}, 2.15, 2.3 and 2.1, Dold's indices of $f$ and of $\hat f|V $ (cf. 5.5) agree and are defined to be the value of $1\in \pi_{stable}^0(B^+)$ under the induced homomorphism of $\hat u$. In other words, we can represent $I(f)$ by the obvious composite map

$$ u:\Sigma^k B^+ @> {\hat u} >> \Sigma^kB^+ \to \Sigma ^k({\{point\}^+})=S^k. \leqno(5.7)$$

\bigskip

 Let us apply the Pontrjagin-Thom procedure (as described above) to $I(f)=[u]$. Clearly $u^{-1}(\{0\})$ is just the fixed point set $C$ of $f$ (cf. 5.3). The trivialization

$$\bar g_B: \nu(C, B\times \Bbb R^k)=\nu_1\oplus \nu_2|C@>\cong >> TF(p_M)\oplus\nu_2|C=C\times \Bbb R^k$$
\bigskip

\noindent is induced by the tangent map of $id-f$. On $\nu_1$ it coincides with $\bar g_B^{\#}$ (cf. 5.4) and on $\nu_2$ it is given by the identity map (since $f$ is constant along each normal ball in the tubular neighborhood $V$ of $M$ in $B \times \Bbb R^k)$. Thus the data $(C\subset B \times \Bbb R^k\to B, \ \bar g_B), \ \  k>>0$, which describe $PT(I(f))$ are just the stabilized coincidence data of $(id,f)$, projected down to $B$. This proves the identity claimed in theorem 1.2. \quad
\qed

{\it Remark 5.8.} A key point in the  previous proof is the fact that Dold's index remains unchanged by the passage $f \rightsquigarrow \hat f$ (cf. 5.5). This parallels closely the stabilizing transition $\omega^{\#}(id,f) \rightsquigarrow \  \tilde \omega (id,f)$.

\head 6.  $S^1-$bundles over $S^1$  \endhead

In this section we study the example 1.22 of the introduction in some detail. In particular, we prove theorems 1.3 and 1.4.

 Given possible values  $q$ and $r$ of the numerical invariants discussed in Proposition 1.1, let us describe the corresponding map $f:M \to N$ in standard form: for any element $[(t,z)]$ in the domain  (cf. 1.23 or 1.24), the standard map is defined by

$$ f([(t,z)])=\cases 
[(t,e^{2\pi irt}z^q)] &\ \hbox{if}\ N=T;\\
[(t,(-1)^rz^q)] &\ \hbox{if}\ N=K.\\
\endcases
\leqno(6.1)
$$

\bigskip

Using  the linear structures on    the universal covering spaces of $T$ and $K$ we see that every map over $S^1$ can be deformed  over $S^1$ into standard form. This proves Proposition 1.1.

 Next we calculate the group $\Omega_1(M;\varphi)$ in which the (weakened) coincidence invariant $\omega_B(f_1,f_2)$ of maps $f_1,f_2:M \to N$ over $S^1$ lies. Here $\varphi$ is trivial if $M=N$; $\varphi$ is the pullback $p^*_M(\lambda)$ of the nontrivial line  bundle over $B=S^1$ if $M \ne N$ (cf. 1.15).

 From \cite{7}, theorem 9.3, we obtain the exact sequence

$$ 0\to \Omega_1^{fr} @> {\delta}>> \Omega_1(M;\varphi)@> {\gamma}>> \bar\Omega_1(M;\varphi) \to 0
\leqno (6.2)$$

\bigskip

\noindent where $\gamma$ forgets about stable vector bundle isomorphisms and retains only  the corresponding orientation information. If $M=N$, then $\gamma$ maps the classical framed bordism group 
$\Omega_1^{fr}(M)$ to the oriented bordism group $\Omega_ 1(M)\cong H_1(M;\Bbb Z)$, and the obvious forgetful homomorphism $\Omega_1^{fr}(M)\to \Omega_1^{fr}$ yields a splitting of  6.2. If $M\ne N$ then a splitting can be extracted from the exact Gysin sequence

$$\Omega_1^{fr}(M) @> d >>\Omega_1^{fr}(\tilde M)@> {proj_*}>>\Omega_1(M;\varphi) \to 0 \leqno (6.3)$$

\bigskip

\noindent where $\tilde M$ is the double cover (or $S^0$-bundle) corresponding to the line bundle 
$\lambda_M:=p^*_M(\lambda)$ over $M$, $d$ takes double coverings and proj denotes the obvious projection. (This is essentially the exact sequence of the pair $(\lambda_M, \lambda_M-s_0(M)$) and uses the Thom isomorphism

$$\Omega_i(\lambda_M,\lambda_M-s_0(M);-\lambda_M)\cong\Omega_{i-1}^{fr}(M)$$

\bigskip

\noindent obtained by intersecting transversely with the zero section $s_0$ of $\lambda_M$).

 Recall that any connected closed smooth $1$-manifold $S$ can carry two distinct stable framings:

(i) the invariant framing obtained from a {\it nonstable} parallelization $TS\cong S\times \Bbb R$ (which is essentially invariant under rotations along the circle $S\cong S^1$); and 

(ii) the boundary framing induced from a disk $D$ which bounds $S=\partial D.$

\noindent  The corresponding bordism classes are $1$ and $0$, resp., in $\Omega_1^{fr}\cong \Bbb Z_2$.

 Now we can describe the direct sum decomposition of $\Omega_1(M;\varphi)$ in theorem 1.4. The projection 
to the first (and the second, resp.) component group is obtained via intersecting circles in $M$ with the fibre $F_M$ 
(and with the section $s_{-1}(B)$ at -1, resp.); the three direct summands are generated by the circles $s_{+1}(B)$ and 
$F_M$ (both with the boundary framing) and by

$$\delta(1) \ := \ [(invariantly \ framed \ S^1, \  constant \  map)]. \leqno (6.4)$$

\bigskip

 In order  to compute the summands  of $\omega_B(f_1,f_2)$ (corresponding to this  decomposition of $\Omega_1(M;\varphi)$)  we may assume that $f_1$, $f_2$ are in standard form (cf. 6.1). Then the pairs $(f_1,f_2)$ and $(f:=f_1\circ f_2^{-1}$, $f_0:=f_2\circ f_2^{-1}=s_{+1}\circ p_M)$ have
the same coincidence locus $C$ which consists of ``parallel" circles in $M$. (Here we use fibrewise complex multiplication of standard maps; it is compatible with the gluing diffeomorphisms of $T$ and $K$, cf. 1.23 and 1.24). The transverse intersections of $C$ with $F_M$ and $s_{-1}(B)$ determine $q$ and $r$ (as indicated in the first two columns concerning 
$\omega_B(f_1,f_2)$ in the table of theorem  1.4; the correction terms 1 and $\rho_2(q)$ result from the fact that the sections $s_{+1}$ and  
$s_{-1}$ have each a self-intersection in $K$).

 Furthermore each circle $S$ in the coincidence locus $C$ is invariantly framed and hence contributes nontrivially to the third component of $\omega_B(f_1,f_2)$; it constitutes a full Nielsen class which therefore must be essential (see also the following proof). This establishes theorem 1.4.



{\bf Proof of Theorem 1.3.} If $N=S^1\times S^1,$ we are in the special case of a product fibration (c.f. 1.17), and  our coincidence theory of maps $f_1,f_2$ over $B$ reduces to the classical coincidence theory of their projections
 $ f_1',  f_2'$ to the fibre $S^1$. But this situation has been thoroughly discussed in \cite{9}, theorem 1.13 and section 6, where even the fibre homotopy type  of $E( f_1',  f_2')$ over $M$ is described. In particular, the Reidemeister number is just the cardinality of the cokernel of the induced homomorphism
$$ f_{1*}'- f_{2*}':H_1(M, \Bbb Z) \to H_1(S^1,\Bbb Z) \cong \Bbb Z$$
\noindent whose image is generated by the greatest common divisor of $(q(f_1)-q(f_2))$ and $(r(f_1)-r(f_2))$. The Reidemeister number equals $MCC(f_1',f_2')=N(f_1'
,f_2')$ except in the selfcoincidence case $f_1 \sim_B f_2$
when $(f_1,f_2)$ is loose (cf. \cite{9}, 1.13). 

If $N=K$ the only sections (up to homotopy) of $p_N$ are $s_{\epsilon}$, $\epsilon=\pm 1$ (cf. 1.25);
each can be deformed away from itself until it has only one selfintersection point in $K$. Therefore,
if maps  $f_i: M \to K$ over $S^1$ are homotopic to $s_{\epsilon_i}\circ p_M$, $i=1,2$, (e.g. if $M=T$), 
their coincidence data can be represented by a whole fibre (or by $\emptyset$ , resp.) when 
$\epsilon_1=\epsilon_2$ (or $\epsilon_1\ne \epsilon_2$, resp.), and $MCC_B(f_1,f_2)=N_B^{\#}(f_1,f_2)= N_B(f_1,f_2)$ equals 1 (or 0, resp.).

It remains to study the coincidence behaviour of maps $f_1,f_2:K \to K$ over $S^1$ in standard form or, equivalently, of maps $f_1 \circ f_2^{-1}=:f$ and $f_2 \circ f_2^{-1}=s_{+1}\circ p_K=:f_0$ (here we use fibrewise complex  multiplication). In view of the previous paragraph we may assume $q\ne 0$. Then the locus $C(f_1,f_2)$ consists of ''horizontal'' circles which are  ''parallel'' to the  sections $s_{\pm 1}$ and intersect each fibre $S^1$ in 
$\eta, \eta z_1,....,\eta z_1^{|q|-1}$   where $z_1=e^{2\pi i/|q|}$ and $\eta=e^{\pi ir/|q|}$  for $r=0,1$.
                                                                                                                                                  
Given $0\leq k < k' < |q|$, we need to know when the coincidence points $\eta z_1^{k},\eta z_1^{k'}$ are  Nielsen equivalent over $B$. This happens precisely  if there is a path 
 $c=l \cdot \hat c$ from $\eta z_1^{k}$ to $\eta z_1^{k'}$ in $K$ 
(consisting of a loop $l$ at $ \eta z_1^{k}$, followed by a path $\hat c$ in the fibre) such that 
$f\circ c=(f\circ l)\cdot (f\circ \hat c)$ is homotopic  to $f_0 \circ c \sim f_0\circ l$ keeping end points fixed. In other words,
the loop $f\circ \hat c$ which winds $k'-k+jq,$ $j \in \Bbb Z$, times around the fibre is Reidemeister equivalent to the trivial loop. Since $\pi_2(S^1)=0$ (cf. 3.8) this means that  $k'$ is equal to $k$ or 
to  $-k+((-1)^r-1)/2$ mod $|q|$, or, equivalently, that $\eta z_1^{k}$ and $\eta z_1^{k'}$
lie in the same coincidence circle (due to the glueing reflection of $K$). Thus each Reidemeister class corresponds to an essential Nielsen class which consist of a single ''horizontal`` circle with winding number 
$\pm 1$ or $ \pm 2$ with respect to the base $S^1$.

Recall that the Reidemeister numbers were computed in Proposition 3.1.

\quad\qed\

{\it Remark 6.5.} It is intriguing to compare the roles of the involution $\iota$ (in the proof of proposition 3.1) on the one side and of complex conjugation (in the proof above) on the other side.

\enddocument

\newpage

Let   $f_1,f_2:M \to K$ be a pair of maps over $S^1$

Let us   compute  the Reidemeister number over $S^1$  of a pair of maps $f_1,f_2:M \to K$ over $S^1$ where $M$ is either 
$T$ or $K$. The case where $M=T$ is simpler than the case where $M=K$  and we will leave to the reader (compare \cite{3}).
Suppose that  $f_{1\#}(a_K)=a_K^{q_1}$ , 
$f_{1\#}(c_K)=a_K^{r_1}c_K$ and  $f_{2\#}(a_K)=a_K^{q_2}$, $f_{2\#}(c_K)=a_K^{r_2}c_K$. Recall from the introduction that $q=q_1-q_2$ and 
$r=r_1-r_2$.

\proclaim{Proposition 3.1} The Reidemeister classes $R_B(f_1,f_2)$ is the set of equivalent classes of elements of the integers $\Bbb Z$  given as follows:  the elements  $m,n\in \Bbb Z$ are in the same Reidemeister clases over $S^1$  if and only if either $m-n\equiv 0$ mod $q$ or  $m+n\equiv r$ mod $q$.    

 i) If $q=0$ then the Reidemeister class over $S^1$ 
which contains an integer  $l$ is the subset $\{l,l-r\}$.  So the number of Reidemeister classes is always infinite.

 ii) If $q\ne 0$ then we have that the number of Reidemester classes over $S^1$ is finite and:

$$\#R_B(f_1,f_2)= \cases 
|q|/2 &\ \hbox{if}\ q \ even, \ r=1\\
[\ |q|/2\ ]+1 &\ \hbox{else}.
\endcases
$$
\endproclaim

\demo{Proof}  A Reidemeister  
 class over $S^1$ of an element $\alpha=a_K^m\in \pi_1(S^1)$ is given by the set of elements of the form 
$$f_{2\#}(a_K^tc_K^v)a_K^m(f_{1\#}(a_K^tc_K^v))^{-1}=a_K^{q_2t}(a_K^{r_2}c_K)^va_K^m(a_K^{q_1t}(a_K^{r_1}c_K)^v)^{-1}=$$
\noindent$=a_K^{qt+ r(1-(-1)^v)/2+(-1)^{ v}m},$ where $t,v \in \Bbb Z$. 
If $q=0$ then each Reidemeister class over $S^1$ contains 2 elements and we have 
 the number of Reidemeister classes   infinite. If $q\ne 0$ then we have $q$  usual Reidemeister classes on the fibre  and a Reidemeister class over $S^1$ is given as the union of 2 such Reidemeister  classes (the   Reidemeister class which contains $m$ and the  Reidemeister class which contains $r-m$). These two  Reidemeister classes are the same presicely in the following  cases: a)  If $r$ is even and $q$ is even then the Reidemeister class which contains  $r/2$ (the one  which contains $q/2+r/2$) coincide with the Reidemeister class over $S^1$ which contains
$q/2$(the one which contains $q/2+r/2$). b)  If $r$ is even and $q$ is odd then the    Reidemeister classes which contains $r/2$ coincide with the Reidemeister class over $S^1$ which contains
$r/2$. c)   If $r$  and $q$ are odd then the    Reidemeister class which contains $r(q+1)/2$ coincide with the Reidemeister class over $S^1$ which contains $r(q+1)/2$. d) If $r$ is odd  and $q$ is even   then a   Reidemeister class over 
$S^1$ is always the union of two distinct  Reidemeister classes. These cases follow from the analyzes of the equation
 $2m\equiv r$ mod $q$. o the result follows.

\quad\qed

\enddemo

\enddocument